\input amstex
\mag=\magstep1
\documentstyle{amsppt}
\nologo
\NoBlackBoxes
\rightheadtext\nofrills{\eightpoint Addition law attached to a stratification}
\leftheadtext\nofrills{\eightpoint V.Z.Enol'skii, S.Matsutani, and Y \^Onishi}
\hcorrection{8pt}
\vcorrection{-30pt}
\define\dint{\dsize\int}
\def\fp{\flushpar}
\define\inbox#1{$\boxed{\text{#1}}$}
\define\underbarl#1{\lower 1.4pt \hbox{\underbar{\raise 1.4pt \hbox{#1}}}}
\define\ce#1{\lceil#1\rceil} 
\define\tp#1{\negthinspace\ ^t#1}
\define\dg#1{(d^{\circ}\geqq#1)}
\define\Dg#1#2{(d^{\circ}(#1)\geqq#2)}

\define\jrac#1#2{\dfrac{\lower 2pt\hbox{$#1$}}{\raise 2pt\hbox{$#2$}}}
\define\sz#1{^{\raise 1pt\hbox{$\sssize #1$}}\hskip -1pt}
\define\lr#1{^{\sssize\left(#1\right)}}
\define\br#1{^{\sssize\left[#1\right]}}
\define\do#1{^{\sssize\left<#1\right>}}

\define\twolower#1{\hbox{\raise -2pt \hbox{$#1$}}}

\define\hookdownarrow{\kern 0pt
 \hbox{\vbox{\offinterlineskip \kern 0pt
  \hbox{$\cap$}\kern 0pt
  \hbox{\hskip 3.5pt$\downarrow$}\kern 0pt}}
}  

\def\maru#1{{\ooalign{\hfil#1\/\hfil\crcr\raise.167ex\hbox{\mathhexbox20D}}}} 
\font\thirteenptmbit=cmmib9 scaled \magstep2
\font\twelveptmbit=cmmib8 scaled \magstep 2
\font\tenptmbit=cmmib10 
\font\eightptmbit=cmmib8
\font\sixptmbit=cmmib6
\define\thirteenptbk#1{\text{\thirteenptmbit{#1}}}
\define\twelveptbk#1{\text{\twelveptmbit{#1}}}
\define\bk#1{\text{\tenptmbit #1}}
\define\eightptbk#1{\text{\eightptmbit{#1}}}
\define\sixptbk#1{\text{\sixptmbit{#1}}}



\font\thirteenptbf=cmbx8 scaled \magstep3



%
\redefine\jtag#1{\tag"$B!J(B#1$B!K(B"}
\def\boxit#1{\vbox{\hrule\hbox{\vrule\kern3pt
    \vbox to 43pt{\hsize 182pt\kern3pt#1\eject\kern3pt\vfill}
    \kern1pt\vrule}\hrule}} 
\def\oVrule{\vrule width .5pt}
\def\oHrule{\hrule height .5pt}
\def\lbox#1#2#3{\kern 0pt 
\dimen1=#1pt \dimen2=#2pt%
\advance\dimen1 by -1.0pt
\dimen3=\dimen1%
\advance\dimen3 by -16pt
\advance\dimen2 by -1.0pt
\dimen4=\dimen2%
\advance\dimen4 by -20pt
\vbox to #1pt
{
\hsize #2pt
\oHrule
\hbox to #2pt
{
\vsize \dimen1
\oVrule
\vbox to \dimen1
{
\hsize \dimen2
\vskip 8pt
\hbox to \dimen2
{
\vsize \dimen3
\hskip 10pt
%
%
\vbox to \dimen3{\hsize \dimen4\fp #3\vfil}%
\hfil\hskip 10pt
}%
\vskip 8pt
}%
\oVrule%
}%
\oHrule%
}\kern 0pt
}%
\def\cbox#1#2#3{\kern 0pt 
\dimen1=#1pt \dimen2=#2pt%
\advance\dimen1 by -1.0pt
\dimen3=\dimen1%
\advance\dimen3 by -3pt
\advance\dimen2 by -1.0pt
\dimen4=\dimen2%
\advance\dimen4 by -3pt
\vbox to #1pt%
{%
\hsize #2pt
\oHrule%
\hbox to #2pt%
{%
\vsize \dimen1%
\oVrule%
\vbox to \dimen1%
{%
\hsize \dimen2%
\vskip 1.5pt
\hbox to \dimen2%
{%
\vsize \dimen3%
\hskip 1.5pt\hfil
\vbox to \dimen3{\hsize \dimen4\vfil\hbox{#3}\vfil}%
\hfil\hskip 1.5pt
}%
\vskip 1.5pt
}%
\oVrule%
}%
\oHrule%
}\kern 0pt
}%
\def\cboxit#1#2#3{$\hbox{\lower 2.5pt \hbox{\cbox{#1}{#2}{#3}}}$}
\define\Sym{\hbox{\rm Sym}}

\document
\baselineskip=12pt 

\topmatter

\title\nofrills{\thirteenptbf  The Addition Law Attached to a Stratification of a 
  Hyperelliptic Jacobian Variety}
\endtitle
\author
  Victor Enolskii, Shigeki Matsutani,Yoshihiro \^Onishi
\endauthor
\address 
  Victor Enolskii\newline
  Department of Mathematics and Statistics,\newline
  Concordia University, 7141 Sherbrooke West, Montreal H4B 1R6,\newline
  Quebec, Canada \newline
  {\rm Email address} : {\tt vze\@ma.hw.ac.uk}
\endaddress
\address 
  Shigeki Matsutani\newline
  8-21-1 Higashi-Linkan, Sagamihara, 228-0811, JAPAN \newline
  {\rm Email address} : {\tt rxb01142\@nifty.com}
\endaddress
\address
  Yoshihiro \^Onishi\newline
  Faculty of Humanities and Social Sciences\newline
  Iwate University, Ueda 3-18-34, Morioka 020-8550, JAPAN \newline
  {\rm Email address} : {\tt onishi\@iwate-u.ac.jp}
\endaddress

\subjclass 
  Primary 14H05, 14K12; Secondary   14H45, 14H51
\endsubjclass
\keywords 
  Schottky-Klein formulae, hyperelliptic sigma functions,
  a subvariety in a Jacobian 
\endkeywords

\abstract
This article shows explicit relation between fractional expressions of Schottky-Klein type 
for hyperelliptic $\sigma$-function 
and a product of differences of the algebraic coordinates on each stratum of natural stratification  
in a hyperelliptic Jacobian. 
\endabstract
\endtopmatter


\document

\head  1. Introduction \endhead

In this paper we shall consider the addition law on the Jacobian 
for any hyperelliptic curve.  
In the theory of Abelian functions, 
addition laws are definitive for a class of corresponding functions. 

To describe our investigation, we start from genus one case, 
and recall that Weierstrass showed that principal relations 
in the theory of elliptic functions can be derived from 
the well-known addition formula
  $$
  \frac{\sigma(v+u)\sigma(v-u)}{\sigma(v)^2\sigma(u)^2}
  = \wp(v) - \wp(u),
  \tag{1.1}
  $$
where $\wp(u)=-(d^2\hskip -2pt/du^2)\log\sigma(u)$.
%

If we consider an elliptic curve  $C_1$  defined by  $y^2=x^3+\lambda_2x^2+\lambda_1x+\lambda_0$  
with unique point  $\infty$  at infinity.  
Let  $(x(u),y(u))$  be the inverse function of 
  $$
  (x,y)\mapsto  u=\int_{\infty}^{(x,y)}\frac{dx}{2y} \hskip 10pt \text{modulo the periods}.
  $$
Then  $x(u)$  equals to the function  $\wp(u)$  attached to  $C_1$  
up to an additive constant.  Hence we have 
  $$
  \frac{\sigma(v+u)\sigma(v-u)}{\sigma(v)^2\sigma(u)^2}
  = x(v) - x(u).
  \tag{1.2}
  $$
\par
We recall two or three kinds of generalizations of (1.1) or (1.2). 
While for the sake of simplicity we restrict here materials only about genus two curve 
$C_2$  defined by  $y^2=x^5+\lambda_4x^4+\lambda_3x^3+\lambda_2x^2+\lambda_1x^1+\lambda_0$, 
in many cases such generalizations are proved for any hyperelliptic curve.    
Before describing our generalizations, 
we note that there is a nice generalization of the Weierstrass elliptic  $\sigma$-function
(see 3.1 below for the definition),  
which is a theta function on  $\Bbb C^2$ (now we are assuming the genus is two).  
Indeed if we define  $\wp_{ij}(u)=-(\partial^2\hskip -2pt/\partial u_i\partial u_j)\log\sigma(u)$, 
then we have a classical formula
\footnote{ Generalization of this formula to the case $g=3$ was given by Baker 
\cite{{\bf B2}}. Recently in \cite{{\bf BEL}} the addition law 
for all genera $g>1$ of hyperelliptic curve 
was written in terms of pffafian built in Kleinian $\wp$-functions.}
  $$
  \frac{\sigma(u+v)\sigma(u-v)}{\sigma(u)^2\sigma(v)^2}=\wp_{11}(v)
  -\wp_{11}(u)-\wp_{12}(u)\wp_{22}(v)+\wp_{12}(v)\wp_{22}(u),
  \tag{1.3}
  $$
which is a natural generalization of (1.1) based on the idea of Jacobi.  
This formula is expressed in terms of the algebraic coordinates as 
  $$
  \aligned
  \frac{\sigma(u+v)\sigma(u-v)}{\sigma(u)^2\sigma(v)^2}
  &=\frac{f(x_1,x_2)-2 y_1y_2}{(x_1-x_2)^2}
  -\frac{f(z_1,z_2)-2 w_1w_2}{(w_1-w_2)^2}\\
  &+ x_1x_2(z_1+z_2)
  - z_1z_2(x_1+x_2),
  \endaligned
  \tag{1.4}
  $$
where  $u=\big(\int_{\infty}^{(x_1,y_1)}+\int_{\infty}^{(x_2,y_2)}\big)(\frac{dx}{2y}, \frac{xdx}{2y})$,  
       $v=\big(\int_{\infty}^{(z_1,w_1)}+\int_{\infty}^{(z_2,w_2)}\big)(\frac{dx}{2y}, \frac{xdx}{2y})$, and 
$f(x, z)$  is a rather complicated polynomial of  $x$  and  $z$  
such that  $f(x,z)=f(z,x)$ (see \cite{{\bf B1}, p.211} for this).

On the other hand, the following is another generalization of (1.2) given by the third author.
Let  $u$  and  $v$  vary on the canonical universal Abelian covering of 
the curve  $C_2$  presented in  $\Bbb C^2$ (this is no other than  $\kappa^{-1}(\Theta\br{1})$  
by notation in Section 2 below). 
and  $x(u)$  be a function such that  $u=\int_{\infty}^{(x(u),y(u))}(\frac{dx}{2y},\frac{xdx}{2y})$.   
Then we have
  $$
  -\frac{\sigma(v+u)\sigma(v-u)}{\sigma_2(v)^2\sigma_2(u)^2}
  =x(v) - x(u),
  \tag{1.5}
  $$
where  $\sigma_2(u)=\partial\sigma(u)/\partial u_2$


Now, it is natural to want a unified understanding of whole the formulae (1.3), (1.4), and (1.5) above.  
There are several hints. 
The first hint would be the fact that the right hand side of (1.3) has 
a determinant expression of \cite{{\bf EEP}},  
which is given by using the genus two case of results in \cite{{\bf \^O}}. 
The second hint is the following formula: 
  $$
  \multline
  -\frac{\sigma(u\lr{1}+u\lr{2}+v) \sigma(u\lr{1}+u\lr{2}-v)}
  {\sigma(u\lr{1}+u\lr{2})^2 \sigma_2(v)^2} \\
   =x(v)^2-\wp_{22}(u\lr{1}+u\lr{2})x(v)-\wp_{12}(u\lr{1}+u\lr{2})
  \endmultline
  \tag{1.6}
  $$
for  $u\lr{1}$, $u\lr{2}$, and  $v$  varying 
on the canonical universal Abelian covering of  $C_2$  in  $\Bbb C^2$.  
This appeared in \cite{{\bf G}} at the first time, 
and a generalization of this for any hyperellitpic curve was reported in \cite{{\bf BES}, (3.21)},  
without proof. 


The main result (Theorem 4.2 below) of this paper seems to be a unification of (1.5) and (1.6). 
There should exist a unification of all the formulae above, 
and such a formulation might be given in near future.  

\vskip 20pt

\fp{\bf Notations}. The symbol $\big(d^{\circ}(z)\ge n\big)$ denotes 
terms of total degree at least  $n$  with respect to a variable  $z$. \fp

\newpage

\head {2. Staratification of the Jacobian of a Hyperelliptic Curve} \endhead

Throughout this article we deal with a hyperelliptic (or elliptic) curve $C_g$  of genus
$g>0$  given by the affine equation
  $$
  y^2=f(x),
  $$
where we are assuming that  $f(x)$  is of the form
  $$
  \aligned
  f(x)&= x^{2g+1} + \lambda_{2g} x^{2g}+\cdots +\lambda_2 x^2 +\lambda_1 x+\lambda_0  \\
  \endaligned
  $$
with  $\lambda_j$'s  being complex numbers.
Then a canonical basis of the space of the differentials of the first kind on  $C_g$  is given by
  $$
  \omega_1:=\frac{       dx}{2y}, \quad
  \omega_2:=\frac{x      dx}{2y}, \quad \cdots, \quad
  \omega_g:=\frac{x^{g-1}dx}{2y}
  $$
with the algebraic coordinate  $(x, y)$  of  $C_g$,
Let  $\alpha_j$  and  $\beta_j$ ($j=1$, $\cdots$, $g$) be a standard homology basis on  $C_g $.  
Namely, they give
  $$
  \text{\rm H}_1(C_g, \Bbb Z)
  =\bigoplus_{j=1}^g\Bbb Z\alpha_{j}
   \oplus\bigoplus_{j=1}^g\Bbb Z\beta_{j},
  $$
and their intersection products are given by  
$[\alpha_i, \alpha_j]=0$, $[\beta_i, \beta_j]=0$, 
$[\alpha_i, \beta_j]=-[\beta_i, \alpha_j]\delta_{i,j}$.
We denote matrices of the half-periods with respect to the differentials  $\omega_i$  
and the homology basis  $\alpha_j$, $\beta_j$  by
  $$
  \omega' :=\frac{1}{2}\left[\int_{\alpha_j}\omega_i\right],
  \quad
  \omega'':=\frac{1}{2}\left[\int_{\beta_j }\omega_i\right].
  $$
We introduce differentials of the second kind
  $$
   dr_j:=\dfrac{1}{2y}\sum_{k=j}^{2g-j}(k+1-j)
    \lambda_{k+1+j} x^kdx,
     \quad (j=1,\ldots,g)
  $$
and the half-periods matrix 
  $$
  \eta':=\frac{1}{2}\left[\dint_{\alpha_{j}}dr_{i}\right], \quad
  \eta'':=\frac{1}{2}\left[\dint_{\beta_{j}}dr_{i}\right] 
  $$
of this differentials with respect to  $\alpha_j$  and  $\beta_j$. 
These  $2g$  meromorphic differentials  $u_i$  and  $r_i$ ($i=1$, $\cdots$, $g$) 
are chosen in such the way that 
the half-periods matrices  $\omega'$, $\omega''$, $\eta'$, $\eta''$  satisfy
the generalized Legendre relation
  $$
  {\frak M} \left[\matrix 0   & -1_g \\ 
                          1_g &  0   
                 \endmatrix\right]
  \tp{\frak M}=\frac{\sqrt{-1}\pi}{2}
            \left[\matrix 0   &  -1_g \\
                          1_g &   0
            \endmatrix\right],
  $$
where  ${\frak M}=\left[\matrix \omega' & \omega'' \\ 
                                \eta'   & \eta''   
                  \endmatrix\right]$.
Let  $\Lambda=2\big(\Bbb Z^g\omega'\oplus \Bbb Z^g\omega''\big)$. 
Then  $\Lambda$  is a lattice in  $\Bbb C^g$  and 
the Jacobi variety  $\Cal J(C_g)$  of  $C_g$  is given by
  $$
  \Cal J(C_g):= \Bbb C^g/\Lambda. 
  $$
We use the modulus  $\Bbb T:={\omega'}^{-1}\omega''$  to define 
the $\sigma$-function of  $C_g$  later.  

For  $k=0$, $\cdots$, $g$, 
the {\it Abel map} $\phi_k$  from the  $k$-th symmetric product  $\Sym^k(C_g)$  
of the curve  $C_g$  to  $\Cal J(C_g)$  is the map
  $$
  \phi_k:\Sym^k(C_g) \rightarrow  \Cal J(C_g) \  \ \text{given by} \ \  \ 
  (Q_1,\ldots,Q_k)\mapsto 
  \sum_{i=1}^k \int_\infty^{Q_i}(\omega_1, \cdots, \omega_g) \mod \Lambda.
  $$
We denote the modulo $\Lambda$  map by  $\kappa$:
  $$
  \kappa:\bold C^g \rightarrow \Cal J(C_g)=\Bbb C^g/\Lambda.
  $$
We denote  $\Theta\br{k}$  the image  $\phi_k(\Sym^k(C_g))$  of the Abel map  $\phi_k$  above. 
Now we have the following stratification:
  $$
  \{O\}=\Theta\br{0}\subset\Theta\br{1} \subset \Theta\br{2} \subset \cdots \subset
  \Theta\br{g-1} \subset \Theta\br{g}=\Cal J(C_g), 
  $$
where  $O$  is the origin of  $\Cal J(C_g)$. 
It is known that each  $\Theta\br{k}$  is a subvariety of  $\Cal J(C_g)$.  
We shall refer each subvariety  $\Theta\br{k}$  by $k$-th stratum of  $\Cal J(C_g)$.

The following Lemma is shown by a straightforward calculation of the Abelian integral 
  $$
  u_i=\int_{\infty}^{(x,y)}\frac{x^{i-1}dx}{2y}
  $$
by using power series expansion and integral by its terms. 

\proclaim{\bf Lemma 2.1}
Let  $u=(u_1, \cdots, u_g)\in\kappa^{-1}(\Theta\br{1})$.  
We denote by  $\big(x(u),y(u)\big)\in C_g$  the algebraic coordinate such that 
whose image by the Abel map  $\phi_k$  is  $u$  modulo  $\Lambda$.  
Then we have the follwoing properties. \fp
{\rm (1)} The variable  $u_g$  is a local parameter around  $(0,\cdots,0)$  on  $\kappa^{-1}(\Theta\br{1})$,
and so that  $u_1$, $\cdots$, $u_g$  are functions of  $u_g$  on  $\kappa^{-1}(\Theta\br{1})$  
near  $(0,\cdots,0)$.\fp
{\rm (2)} The function  $x(u)$  has the following Laurent expansions
around $(0,0,\cdots,0)$  on  $\kappa^{-1}(\Theta\br{1})$ {\rm :}
  $$
  \align
  x(u)&={u_g}^{-2} +  \big(d^{\circ}(u_g) \ge 0\big), \\
  y(u)&={u_g}^{-2g-1} + \big(d^{\circ}(u_g) \ge -2g+1\big).
  \endalign
  $$
\endproclaim

For a proof of the above, we refer the reader to Lemmas 3.8 and 3.9 in \cite{{\bf \^O}}, for instance. \qed

\head 3. The Sigma Function and Its Derivatives \endhead

In this section, we will introduce the hyperelliptic $\theta$-function
and $\sigma$-function. The later one is a natural generalization of the
Weierstrass $\sigma$-function.

Let $a$  and  $b$ be two vectors in $\Bbb R^g$. 
We recall the theta function  with respect to 
the lattice of periods generated by  $1_g$  and  $\Bbb T={\omega'}^{-1}\omega''$  
with characteristic  $\tp{[a\ b]}$, which is a function of  $z\in\Bbb C^g$  defined by
  $$
  \theta\negthinspace\left[\matrix a \\ b \endmatrix\right]\hskip -2pt (z)
    =\theta\negthinspace\left[\matrix a \\ b \endmatrix\right]\hskip -2pt (z; \Bbb T)
    =\sum_{n \in \Bbb Z^g} 
      \exp \bigg[2\pi\sqrt{-1}\Big\{\frac 12 \tp{(n+a)}\Bbb T(n+a) + \tp{(n+a)}(z+b)\Big\}\bigg],
  $$
as usual. 
Let  $\delta' =\tp{\big[\frac12,   \frac12,     \cdots, \frac12\big]}$  
and  $\delta''=\tp{\big[\frac{g}2, \frac{g-1}2, \cdots, \frac12\big]}$.  
Then the half-period vector  $\delta'\omega'+\delta''\omega''$  is no other than 
so-called Riemann constant for  $\Cal J(C_g)$.  

\definition{\bf Definition 3.1} \ 
The  $\sigma$-function (see, for example, \cite{{\bf B1}, p.336}) 
is given by 
  $$
  \sigma(u)
  =\gamma_0 \exp\Big(\hskip -3pt -\frac{1}{2}\tp{u}\eta'{{\omega}'}^{-1}u \Big)
  \vartheta\negthinspace
  \left[\matrix\delta'' \\ \delta' \endmatrix\right]
  \Big(\frac{1}{2}{{\omega}'}^{-1}u ;\Bbb T\Big), 
  $$
where  $\gamma_0$  is a certain non-zero constant depending of  $C_g$,
which is explained in \cite{{\bf BEL} p.32} and \cite{{\bf \^O}, Lemma 4.2}.
We regards the domain  $\Bbb C^g$  where this function is defined to be  $\kappa^{-1}(\Cal J(C_g))$.  
\enddefinition

Following the paper \cite{{\bf \^O}}, we introduce multi-indices  $\natural^n$  
and their associated derivatives  $\sigma_{\natural^n}(u)$  of  $\sigma(u)$  as follows:

\definition{\bf Definition 3.2}
We define that 
  $$
  \natural^n=\cases
  \{n+1, n+3, \cdots, g-1\} & \text{if}\  g-n\equiv 0 \mod\; 2,\\
  \{n+1, n+3, \cdots, g  \} & \text{if}\  g-n\equiv 1 \mod\; 2.
 \endcases
  $$
By using this notation we have partial derivatives of  $\sigma(u)$  associated these multi-indeces, namely,  
  $$
  \sigma_{\natural^n}(u)=\bigg(\prod_{i\in \natural^n}\frac{\partial}{\partial u_i}\bigg)\sigma(u).
  $$
Moreover we denote  $\sharp:=\natural^1$  and  $\flat:=\natural^2$, 
so that  $\sigma_{\sharp}(u)=\sigma_{\natural^1}(u)$  and  $\sigma_{\flat}(u)=\sigma_{\natural^2}(u)$.  
\enddefinition

Several examples of  $\sigma_{\natural^n}(u)$  are given in the following table
\footnote{One can see that numbers appearing in those multy-index  $\natural^n$ 
are naturally related to the Weierstrass gap sequence  $1$, $3$, $5$, $\cdots$, $2g-1$  
at the Weierstrass point  $\infty$  at infinity.
}.
\vskip 5pt
\fp
  \centerline{
  \vbox{\offinterlineskip
    \baselineskip =5pt
    \tabskip = 1em
    \halign{&\hfil#\hfil \cr
     \noalign{\vskip 2pt}
     \noalign{\hrule height0.8pt}
        & genus & \hfil\strut{\vrule depth 4pt}\hfil  & $\sigma_\sharp\equiv \sigma_{\natural^1}$ & $\sigma_\flat\equiv\sigma_{\natural^2}$ & $\sigma_{\natural^3}$  & $\sigma_{\natural^4}$ & $\sigma_{\natural^5}$ &$\sigma_{\natural^6}$&$\sigma_{\natural^7}$ &$\sigma_{\natural^8}$ &$\cdots$\cr
      \noalign{\hrule height0.3pt}
        &$1$    & \strut\vrule  & $\sigma$ & $\sigma$ & $\sigma$ & $\sigma$
                & $\sigma$      & $\sigma$ & $\sigma$ & $\sigma$ &$\cdots$\cr
        &$2$ & \strut\vrule  & $\sigma_2$ & $\sigma$& $\sigma$ & $\sigma$
                         & $\sigma$ & $\sigma$& $\sigma$ & $\sigma$&$\cdots$\cr
&$3$ & \strut\vrule  & $\sigma_2$ & $\sigma_3$& $\sigma$ & $\sigma$
                         & $\sigma$ & $\sigma$& $\sigma$ & $\sigma$&$\cdots$\cr
&$4$ & \strut\vrule  & $\sigma_{24}$ & $\sigma_3$& $\sigma_4$ & $\sigma$
                         & $\sigma$ & $\sigma$& $\sigma$ & $\sigma$&$\cdots$\cr
&$5$ & \strut\vrule  &$\sigma_{24}$& $\sigma_{35}$& $\sigma_4$ & $\sigma_5$
                         & $\sigma$ & $\sigma$& $\sigma$ & $\sigma$&$\cdots$\cr
&$6$ & \strut\vrule  &$\sigma_{246}$&$\sigma_{35}$& $\sigma_{46}$&
                         $\sigma_5$& $\sigma_6$ & $\sigma$& $\sigma$
                         & $\sigma$&$\cdots$\cr
&$7$ & \strut\vrule  &$\sigma_{246}$&$\sigma_{357}$& $\sigma_{46}$&
                         $\sigma_{57}$& $\sigma_6$ & $\sigma_7$& $\sigma$
                         & $\sigma$&$\cdots$\cr
&$8$ & \strut\vrule  &$\sigma_{2468}$&$\sigma_{357}$& $\sigma_{468}$&
                         $\sigma_{57}$& $\sigma_{68}$ &$\sigma_7$& $\sigma_8$
                         & $\sigma$&$\cdots$\cr
&$\vdots$ & \hfil\strut {\vrule depth 5pt}\hfil  & $\vdots$ & $\vdots$ & $\vdots$ & $\vdots$ & $\vdots$ & $\vdots$ & $\vdots$ & $\vdots$ & $\ddots$ \cr
	\noalign{\hrule height0.8pt}
}}
}

For  $u\in \kappa^{-1} \Cal J(C_g)$,  we denote  $u'$  and  $u''$  
the unique elements in  $\Bbb R^g$  such that  $u=2(u'\omega'+u''\omega'')$.  
We introduce a  $\Bbb C$-valued  $\Bbb R$-bilinear form  $L(\ ,\ )$  defined by 
  $$
  L(u,v) = \tp{u}(\eta' u' + \eta'' u'')
  $$
for  $u$, $v\in \kappa^{-1} \Cal J(C_g)(=\Bbb C^g)$.  Let 
  $$
  \chi(\ell)=\exp\left[2\pi \sqrt{-1}(\tp{\ell'}\delta''-\tp{\ell''}\delta')-\pi\sqrt{-1}\tp{\ell'}\ell''\right].
  $$
The following facts are essential for our main result.  

 \proclaim{Proposition 3.3} 
 {\rm (1)} For  $u\in\kappa^{-1}(\Theta\br{n})$  and  $\ell\in\Lambda$, we have 
   $$
   \sigma_{\natural^n}(u+\ell) =
   \chi(\ell)\sigma_{\natural^n}(u)
   \exp L(u+\tfrac12 \ell, \ell).
   $$
 {\rm (2)} Let  $n$  be a positive integer  $n\le g$.
 Let  $v$, $u^{(1)}$, $u^{(2)}$, $\cdots$, $u^{(n)}$
 be elements in  $\kappa^{-1}(\Theta^{[1]})$. 
 If  $u\lr{1}+\cdots+u\lr{n}\not\in \kappa^{-1}(\Theta\br{n-1})$, 
 then the function  $v\mapsto\sigma_{\natural^{n+1}}(u\lr{1}+\cdots+u\lr{n}+v)$  
 has only zeroes at  $v=(0,\cdots, 0)$  modulo  $\Lambda$  of order  $g-n$  and at  $-u\lr{1}$ 
 modulo $\Lambda$  of order $1$.  
 Around $(0,0,\cdots,0)$ 
 we have the following expansion with respect to  $v_g$ {\rm :}
  $$
  \sigma_{\natural^{n+1}}(u\lr{1}+\cdots+u\lr{n} + v)
  =(-1)^{(g-n)(g-n-1)/2} \sigma_{\natural^{n}}(u){v_g}^{g-n} +
  (d^{\circ}(v_g) \ge g-n+1).
  $$
 \fp
 {\rm (3)} For  $v\in\kappa^{-1}(\Theta\br{1})$, 
  $$
  \sigma_{\natural^{1}}(v)= -(-1)^{g(g-1)/2}{v_g}^{g} + (d^{\circ}(v_g) \ge g+2).
  $$
 \endproclaim
\fp{\it Proof.} \ 
The assertion (1) is proved by Lemma 7.3 in \cite{{\bf\^O}}.   
The assertions (2) and (3) are proved by Proposition 7.5 in \cite{{\bf \^O}}.  
\qed
\vskip 6pt
\fp
{\bf Remark 3.4.}\ 
(1) If  $n=g$, the assertion 3.3(1) is no other than 
the classical relation of  $\sigma(u)$  with respect to the translation by any period  $\ell\in\Lambda$  
(see \cite{{\bf B1}, p.286}). \fp
Namely, we have the same relations for translation by any period  $\ell\in\Lambda$  
for the special partial-derivatives  $\sigma_{\natural^n}(u)$  on  $\kappa^{-1}(\Theta\br{n})$ 
as that of  $\sigma$  itself.  
The essence of the proof in \cite{{\bf \^O}} of this fact is 
that the derivative  $\sigma_{\check{\natural}^n}(u)$  
for any proper subset  $\check{\natural}^n$  of  $\natural^n$   
vanishes on the strutum  $\Theta\br{n}$, 
which is proved by investigating Riemann singularity theorem explicitly.  \fp
(2) We see by considering the unification of 3.3(2) and 3.3(3) above 
that it would be natural to define  $\sigma_{\natural^0}(u)=-1$ (a constant function).  \fp
(3) The statements 3.3(2) and (3) complement the Riemann singularity theorem (\cite{{\bf R}}) 
(see also \cite{{\bf ACGH}, p.226-227}) with pointing the orders of vanishing 
on each stratum  $\Theta\br{n}$ in terms of  $\sigma$  function.
\vskip 0pt

\newpage 


\head 4. Main Result \endhead 

In this Section we start to recall the follwoing formula (Lemma 8.1 in \cite{{\bf \^O}}) 
without proof.

\proclaim{\bf Lemmma 4.1} 
Suppose  $u$  and  $v$  be variables on  $\kappa^{-1}(\Theta\br{1})$. 
Then we have  
  $$
  (-1)^{g-1}\frac{\sigma_{\flat}(v+u)\sigma_{\flat}(v-u)}{\sigma_{\sharp}(v)^2\sigma_{\sharp}(u)^2}
  = x(v)-x(u).
  $$
\endproclaim

The following relation is our main theorem and is an extension of both of (3.21) in \cite{{\bf BES}}
and the formula above. 

\proclaim{\bf Theorem 4.2}
Let  $m$  and  $n$  be positive integers such that  $m+n\leq g+1$.
Let
  $$
  u=\sum_{i=1}^m\int_{\infty}^{(x_i,y_i)}  \hskip -10pt(\omega_1,\cdots,\omega_g) \ \in \kappa^{-1}(\Theta\br{m}), \ \ 
  v=\sum_{j=1}^n\int_{\infty}^{(x'_j,y'_j)}\hskip -10pt(\omega_1,\cdots,\omega_g) \ \in \kappa^{-1}(\Theta\br{n})
  $$
Then the following relation holds {\rm :}
  $$
  \frac{\sigma_{\natural^{m+n}}(u+v) \sigma_{\natural^{m+n}}(u-v)}
  {\sigma_{\natural^m}(u)^2 \sigma_{\natural^n}(v)^2}
  =(-1)^{\delta(g,n)}\prod_{i=1}^m\prod_{j=1}^n (x_i - x'_j),
  $$
where $\delta(g,n)=\frac12n(n-1)+gn-1$.
\endproclaim

\vskip 10pt
\fp
{\it Proof.} \ 
We prove the desired formula by induction with respect to  $m$  and  $n$. 
First of all we suppose that  $2g$  points  $u\lr{1}$, $\cdots$, $u\lr{g}$  
and  $v\lr{1}$, $\cdots$, $v\lr{g}$  are given.  
Then by 3.3 we see that both sides of the desired formula are 
function on  $\Theta\br{1}$  with respect to each variable in  $u\lr{i}$s  and  $v\lr{j}$s. 
We let  $u\br{i}=u\lr{1}+\cdots+u\lr{i}$  
and     $v\br{j}=v\lr{1}+\cdots+v\lr{j}$  
for  $0\leq i \leq g$  and  $0\leq j \leq g$.  
If  $m=n=1$, the assertion is just the Lemma 4.1.   
Therefore, the assertion is proved by reducing
  $$
  \aligned
  &\frac{\sigma_{\natural^{m+n+1}}(u\br{m}+v\br{n+1}) \sigma_{\natural^{m+n+1}}(u\br{m}-v\br{n+1})}
  {\sigma_{\natural^m}(u\br{m})^2 \sigma_{\natural^{n+1}}(u\br{n+1})^2} \\
  &\hskip 140pt =(-1)^{\delta(g,n+1)}\prod_{i=1}^m\prod_{j=1}^{n+1} (x_i - x'_j)
  \endaligned
  \tag{4.3a}
  $$
and
  $$
  \aligned
  &\frac{\sigma_{\natural^{m+n+1}}(u\br{m+1}+v\br{n}) \sigma_{\natural^{m+n+1}}(u\br{m+1}-v\br{n})}
  {\sigma_{\natural^{m+1}}(u\br{m+1})^2 \sigma_{\natural^{n}}(u\br{n})^2} \\
  &\hskip 140pt =(-1)^{\delta(g,n)}\prod_{i=1}^{m+1}\prod_{j=1}^n (x_i - x'_j),
  \endaligned
  \tag{4.3b}
  $$
to the formula 
  $$
  \frac{\sigma_{\natural^{m+n}}(u\br{m}+v\br{n}) \sigma_{\natural^{m+n}}(u\br{m}-v\br{n})}
       {\sigma_{\natural^m}(u\br{m})^2 \sigma_{\natural^n}(v\br{n})^2}
  =(-1)^{\delta(g,n)}\prod_{i=1}^m\prod_{j=1}^n (x_i-x'_j).  
  \tag{4.4}
  $$
We denote the index of sign in 3.3(2) by  $\epsilon(g,n)$, 
that is  $\epsilon(g,n)=(g-n)(g-n-1)/2$.  
Then the left hand side of (4.3a) is
  $$
  \aligned
  &\frac{\sigma_{\natural^{m+n+1}}(u\br{m}+v\br{n+1}) \sigma_{\natural^{m+n+1}}(u\br{m}-v\br{n+1})}
        {\sigma_{\natural^m}(u\br{m})^2 \sigma_{\natural^{n+1}}(v\br{n+1})^2} \\
  &=\Big[\sigma_{\natural^{m+n}}(u\br{m}+v\br{n})
         \big\{(-1)^{\epsilon(g,m+n)}\big( v\lr{n+1}_g\big)^{g-m-n}+\cdots\big\} \\
  &\hskip 50pt
   \times\sigma_{\natural^{m+n}}(u\br{m}-v\br{n})
         \big\{(-1)^{\epsilon(g,m+n)}\big(-v\lr{n+1}_g\big)^{g-m-n}+\cdots\big\}
     \Big] \\
  &\hskip 80pt\bigg/  
        {\sigma_{\natural^m}(u\br{m})^2
         \sigma_{\natural^n}(v\br{n})^2
         \big\{(-1)^{\epsilon(g,n)}\big(v\lr{n+1}_g\big)^{g-n}+\cdots\big\}^2      } \\
  &=\frac{\sigma_{\natural^{m+n}}(u\br{m}+v\br{n}) \sigma_{\natural^{m+n}}(u\br{m}-v\br{n})}
         {\sigma_{\natural^m}(u\br{m})^2 \sigma_{\natural^{n}}(u\br{n})^2}
          \Bigg\{(-1)^{g-m-n}\frac{1}{\big(v\lr{n+1}_g\big)^{2m}}+\cdots\Bigg\}
  \endaligned
  \tag{$4.3\text{a}'$}
  $$
by using the Lemma 2.1. 
The right hand side of (4.3a) is
  $$
  \aligned
  &\hskip -20pt
   (-1)^{\delta(g,n+1)}\prod_{i=1}^m\prod_{j=1}^{n+1}(x_i-x'_j) \\
  &\hskip -10pt =(-1)^{\delta(g,n+1)}
    (-x'_{n+1})^m\prod_{i=1}^m\prod_{j=1}^{n}(x_i-x'_j)+\big(d^{\circ}(x_{n+1})\leq m-1\big) \\
  &\hskip -10pt =
    \frac{(-1)^{\delta(g,n+1)+m}}
         {\big(v\lr{n+1}\big)^{2m}}
    \prod_{i=1}^m\prod_{j=1}^{n}(x_i-x'_j)+\big(d^{\circ}(v\lr{n+1})\geqq -2m+1\big)
  \endaligned
  \tag{$4.3\text{a}''$}
  $$
The index of sign of the last in ($\text{4.3a}''$) is 
  $$
  \align
  \delta(g,n+1)+m
  &=\tfrac12(n+1)n+g(n+1)-1+m \\
  &=\tfrac12n(n-1)+n+gn+g-1+m \\
  &\equiv \tfrac12n(n-1)+gn-1+(g-n-m) \mod 2 \\
  &=\delta(g,n)+(g-n-m). \\
  \endalign
  $$
This is equal to the sum of indeces of sign in (4.4) and one in the last factor in ($\text{4.3a}'$).  
Thus, the leading terms of expansions with respect to  $v\lr{n+1}_g$  of the two sides 
completely coincide.  
Till here, the assumption  $m+n\leq g+1$  is not so essential.   
Now we are going to check the divisors of the two sides 
regarding as functions of  $v\lr{n+1}$  modulo  $\Lambda$.  
Using the assumption  $m+n\leq g+1$, 
we can determine the divisors of two sides exactly by Proposition 3.4.  
For the left hand side of (4.3a), the numerator has zeroes 
at  $v\lr{n+1}=(0,0,\cdots,0)$  of order  $2(g-m-n)$  modulo  $\Lambda$, 
at  $\pm\,u\lr{1}$, $\cdots$, $\pm\,u\lr{m}$, $-v\lr{1}$, $\cdots$, $-v\lr{n}$  
of order $1$  modulo  $\Lambda$. 
The denominator has zeroes 
at  $v\lr{n+1}=(0,0,\cdots,0)$  of order  $2(g-n)$  modulo  $\Lambda$, 
at  $-v\lr{1}$, $\cdots$, $-v\lr{n}$  of order $1$  modulo  $\Lambda$.   
Therefore the left hand side of (4.3a) has only pole
at  $v\lr{n+1}=(0,0,\cdots,0)$  of order  $2m$  modulo  $\Lambda$, 
and has zeroes 
at  $\pm\,u\lr{1}$, $\cdots$, $\pm\,u\lr{m}$  of order $1$  modulo  $\Lambda$.   
These pole and zeroes coincide with 
those of the left hand side including their order, 
because, for  $u$  and  $v\in\kappa^{-1}(\Theta\br{1})$, 
we have  $x(u)-x(v)$  if and only if  $u=\pm v$. 
Thus, we have reduced the formula (4.3a) to the equality (4.4). 
The formula (4.3b) is similarly reduced to (4.4).  
Hence, we have proved the assertion.  
\qed 
\vskip 5pt

Baker (\cite{{\bf B1}} and \cite{{\bf B2}}) defined 
  $$
  \wp_{ij}(u)=-\frac{\partial^2}{\partial u_i\partial u_j}\log\sigma(u)
  $$
for  $0\leq i\leq g$, $0\leq j\leq g$, and $u=(u_1, \cdots, u_g)\in\Bbb C^g$, 
which is natutal generalizations of Weierstrass  $\wp$  function.  
As we mentioned before 4.2 
the following special case of  $(m,n)=(g,1)$  in 4.2 appeared in \cite{{\bf BES}, (3.21)}, 
which was the motivation of this paper.  

\proclaim{\bf Corollary 4.5} {\rm (\cite{{\bf BES}, (3.21)})} \ 
We suppose that  $(x,y)$, $(x_1,y_1)$, $\cdots$, $(x_g,y_g)$  are  $g+1$  points on  $C_g$. 
Let  
  $$
  \align
  F_g(x)&=(x-x_1)(x-x_2) \cdots (x-x_g) \\
  u&=\sum_{i=1}^g\int_{\infty}^{(x_i,y_i)}(\omega_1,\cdots,\omega_g)
  \in \kappa^{-1}(\Theta\br{g})(=\Bbb C^g), \\
  v&=\int_{\infty}^{(x,y)}(\omega_1,\cdots,\omega_g)\in \kappa^{-1}(\Theta\br{1}). 
  \endalign
  $$
Then we have the relation 
  $$
  \align
  \frac{\sigma(u+v) \sigma(u-v)}
  { \sigma(u)^2\sigma_\sharp(v)^2}
  &=(-1)^{\frac12g(g+1)}F_g(x) \\
  &=x^g-\wp_{gg}(u)x^{g-1}-\wp_{g,g-1}(u)x^{g-2}-\cdots-\wp_{g1}(u).
  \endalign
  $$
\endproclaim
\vskip 10pt
\fp
{\it Proof.}\ 
The first equality is obvious from 4.2.  
The second equality follow from the fact
  $$
  \wp_{g,g-k+1}(u)=(-1)^{k+1}\hskip -20pt
  \sum_{1\leq i_1< i_2< \cdots < i_k\leq g}
  \hskip -20pt x_{i_1}x_{i_2}\cdots x_{i_k}
  $$
for  $1\leq k \leq g$, which was given by Baker(see \cite{{\bf B3}}, for example).  
\qed
\vskip 10pt


\head 5. Some Remarks \endhead 

We finally remark on related works and future of our result. 
The polynomial  $F_g$  plays the role of the master
polynomial in the theory of hyperelliptic functions. Bolza enabled
to express the polynomial  $F_g$  in terms of Kleinian $\wp$-functions;
in this context we shall call master polynomial $F_g$ as Bolza polynomial.
Its zeroes give solution to the Jacobi inversion problem. In the $2\times2$
Lax representation of dynamic systems associated with a hyperelliptic
curve it plays the role of $U$ polynomial among Jacobi's  $U$, $V$, $W$-triple (\cite{{\bf Mu}}).  
Vanhaecke studied the properties of $\Theta\br{k}$
using $U$, $V$, $W$-polynomials which are constructed on the basis of the
master polynomial  $F_g$ (\cite{{\bf V2}}). In \cite{{\bf BES}}, the authors applied
aforementioned particular case of the addition theorem
with Bolza polynomial to compute the norm of wave function
to the Schr\"odinger equation with finite-gap potential.

A similar polynomial
$F_k(z)\equiv U(z):=(z-x_1)\cdots(z-x_k)$  over  $\Theta\br{k}$  
which plays essential roles in the studies of structures
of the subvarieties (\cite{{\bf AF}} and \cite{{\bf Ma}}), namely, 
  $$
  \frac{\sigma_{\natural^{k+1}}(u - v) \sigma_{\natural^{k+1}}(u + v)}
  {\sigma_{\natural^1}(v)^2 \sigma_{\natural^k}(u)^2} =(-1)^{g-k} F_{k}(x).
  $$
For the case of $k=1$, it appeared in \cite{{\bf \^O}} as
Frobenius-Stickelberger type relation of higher genus,
which determines an algebraic structure of the curve.

We emphasize that study of $\theta$-divisor attract now more attention.
In particular inversion of higher genera hyperelliptic integrals with respect to
the restriction to $\theta$-divisor was recently used in the problem of
analytic description by means of reduction of the Benney system
of conformal map of a  domain in half-plane to its complement
\cite{{\bf BG1}}, \cite{{\bf BG2}}.  
The same method of inversion of an ultraelliptic integral
was used in \cite{{\bf EPR}} to describe motion of double pendulum.

In was shown in series of publications of Vanhaecke
\cite{{\bf V1}}, \cite{{\bf V2}}, Abenda and Fedorov \cite{{\bf AF}}
and one from authors \cite{{\bf Ma}} and others that $\theta$-divisor can
serve as a carrier of integrability.
Grant \cite{{\bf G}}, Cantor \cite{{\bf C}}  found 
algebraic structures on  $\Theta\br{1}$  which is related to
division polynomials whose zeros determines $n$-time points.  
Recently on from the authors and Eilbeck and Previato 
used Grant-Jorgenson formula (\cite{{\bf G}} and \cite{{\bf J}})  
to derive an analogue of the Frobenius-Stickelberger addition formula
for three variables in the case of genus two hyperelliptic curve (\cite{{\bf EEP}}).

Thus we expect that our main theorem has some effects on these
studies based on the Riemann singularity theorem.

\enddocument
\bye